\documentclass[12pt]{article}
\usepackage{amsmath, amsfonts, amssymb}
\usepackage{multirow}
\def\theequation{\arabic{equation}}
\usepackage{graphicx}
\usepackage{color}
\def\theequation{\thesection.\arabic{equation}}
\def\theequation{\arabic{equation}}

\newsavebox{\savepar}

\newtheorem{remark}{Remark}

\newtheorem{example}{Example}
\newtheorem{definition}{Definition}
\newtheorem{theorem}{Theorem}

\def\appendix{\par
\setcounter{section}{0}
   \setcounter{equation}{0}
    \def\@chapapp{APPENDIX}
    \def\thesection{\Alph{section}}
    \def\theequation{A.\arabic{equation}}
    \def\section{

stepcounter{section}
    \@startsection {section}{A}{0pt}{-3.5ex \@plus -1ex \@minus -.2ex}
                   {0.3ex \@plus.2ex}{\normalsize\sffamily\bfseries} }}

\def\eop{\hbox{{\vrule height7pt width3pt depth0pt}}}

\makeatletter
\newcommand{\least}{\let\CS=\@currsize\renewcommand{\baselinestretch}{1}\tiny\CS}

\newcommand{\oneandahalfspacing}{\let\CS=\@currsize\renewcommand{\baselinestretch}{1.2}\tiny\CS}
\newcommand{\doublespacing}{\let\CS=\@currsize\renewcommand{\baselinestretch}{2.5}\tiny\CS}

\usepackage{pst-plot}
\usepackage{amsmath}
\usepackage{amsfonts}
\usepackage{euscript}
\usepackage{oldgerm}
\usepackage{rotating}
\usepackage{topcapt}
   \renewcommand{\baselinestretch}{1.2}

\makeindex

\begin{document}

\newcommand{\namelistlabel}[1]{\mbox{#1}\hfil}
\newenvironment{namelist}[1]{%
\begin{list}{}
{
\let\makelabel\namelistlabel
\settowidth{\labelwidth}{#1}
\setlength{\leftmargin}{1.1\labelwidth}
}
}{%
\end{list}}

\newcommand{\be}{\begin{equation}}
\newcommand{\ee}{\end{equation}}
\newcommand{\dist}{{\rm\,dist}}
\newcommand{\sspan}{{\rm\,span}}
\newcommand{\re}{{\rm Re\,}}
\newcommand{\im}{{\rm Im\,}}
\newcommand{\sgn}{{\rm sgn\,}}
\newcommand{\beano}{\begin{eqnarray*}}
\newcommand{\eeano}{\end{eqnarray*}}
\newcommand{\bea}{\begin{eqnarray}}
\newcommand{\eea}{\end{eqnarray}}

\newcommand{\ba}{\begin{array}}
\newcommand{\ea}{\end{array}}
\newcommand{\hone}{\mbox{\hspace{1em}}}
\newcommand{\hon}{\mbox{\hspace{1em}}}
\newcommand{\htwo}{\mbox{\hspace{2em}}}
\newcommand{\hthree}{\mbox{\hspace{3em}}}
\newcommand{\hfour}{\mbox{\hspace{4em}}}
\newcommand{\von}{\vskip 1ex}
\newcommand{\vone}{\vskip 2ex}
\newcommand{\vtwo}{\vskip 4ex}
\newcommand{\vthree}{\vskip 6ex}
\newcommand{\vfour}{\vspace*{8ex}}
\newcommand{\norm}{\|\;\;\|}
\newcommand{\integ}[4]{\int_{#1}^{#2}\,{#3}\,d{#4}}
\newcommand{\inp}[2]{\langle {#1} ,\,{#2} \rangle}
\newcommand{\vspan}[1]{{{\rm\,span}\{ #1 \}}}
\newcommand{\R} {{\mathbb{R}}}

\newcommand{\B} {{\mathbb{B}}}
\newcommand{\C} {{\mathbb{C}}}
\newcommand{\N} {{\mathbb{N}}}
\newcommand{\Q} {{\mathbb{Q}}}
\newcommand{\LL} {{\mathbb{L}}}
\newcommand{\Z} {{\mathbb{Z}}}

\newcommand{\BB} {{\mathcal{B}}}
\newcommand{\dm}[1]{ {\displaystyle{#1} } }
\def \stackt{{\stackrel{.}{.\;.}\;\;}}
\def \stackb{{\stackrel{.\;.}{.}\;\;}}
\def \olu{\overline{u}}
\def \olv{\overline{v}}
\def \olx{\overline{x}}
\def \olp{\overline{\partial}}
\def\diag{{\;{\rm diag } \; }}
\thispagestyle{empty}

\bibliographystyle{}

\newcommand{\alert}[1]{\fbox{#1}}


\begin{center}
{\Large \bf  Simpson's Paradox and Collapsibility}
\end{center}

\vthree
\begin{center}
{\bf \large P. Vellaisamy}  \\
{\it Department of Mathematics,
Indian Institute of Technology Bombay, \\ Powai,
Mumbai-400 076, India}\\
\end{center}

\vtwo
\noindent {\bf Abstract}. Simpson's paradox and collapsibility are  two closely related concepts in the context
of data analysis. While the knowledge about the occurrence of Simpson's paradox helps a statistician to draw  correct and meaningful conclusions, the concept of collapsibility deals with dimension-reduction aspects, when
 Simpson's paradox does not occur. We discuss in this paper in some detail the nature and the genesis of Simpson's paradox with respect to well-known examples and also various concepts of collapsiblity. The main aim is to  bring out the close connections between these two phenomena, especially with regard to the analysis of contingency tables,
 regression models and a certain measure of association or a dependence function. There is a vast literature
 on these topics and so we focus only on certain aspects, recent developments and some important results  in the above-mentioned areas.

\vtwo

\vone
\section{ Introduction}
 It is well known that statistics or more precisely statistical techniques play an
 important role in addressing some of the problems of a society, an industry and a country. But,
drawing intelligent and correct decisions form the real-life data
is not straightforward. Indeed, we read/hear different or paradoxical
conclusions in several contexts or situations.
So, a layman gets confused and probably believes in the famous
quote  ``Lies, Damned Lies and Statistics," in that
order. In this paper, we will address one of the well-known
paradoxes due to Simpson. Simpson (1951) discussed
how a simple fact about fractions can lead  to a contradictory
conclusions in a  wide variety  of situations. Though some statisticians (Yule (1903))
were aware of these issues in the beginning of 20th century, it is Simpson who
popularized this paradox, earning his name, through analysis of
  real-life data that arise in  several practical applications.

\vone  Simpson's paradox occurs when an observed
association between two random variables, say $X$ and $Y$, gets
reversed after considering the third variable $W$, called a
covariate or a background variable. The situation of having two contradictory conclusions makes this phenomenon
paradoxical. It is one of the most discussed and studied paradoxes
in the statistics literature. The knowledge and the awareness of its occurrence is of importance for the
statistical analysis of the data. It
arises naturally in several areas which include the analysis of
contingency tables, regression models, measures of association,
survival analysis, etc. In this paper, we will discuss some examples of
Simpson's paradox through some real-life data and then discuss some recent results in the areas mentioned above. A closely
related concept, namely collapsibility, is applied whenever the Simpson's paradox does not occur. We present a
survey of some recent results, its relation to Simpson's paradox and discuss  possible directions for future work.

\section{Simpson's Paradox}
\subsection{Real-life situations}

 We start with two examples discussed in the literature.

\begin{example} {\em  (Graduate Admission Data in UC Berkeley).
Let $A\in\{Y, N\}$ denote the admission, $X\in\{M,F\}$ denote the
sex and $D\in\{H,G\}$ denote the department (H=History, G =
Geography). The following data roughly represents the graduate admission
in the two departments $H$ and $G$ (Bickel {\it et. al} (1975),
discrimination suit against UCB ).
\begin{center}
\begin{tabular}{|c|c|lcr|}
\hline
& & & D& \\
\hline
 A&X&H&&G\\
 \hline
\multirow{2}{*}{Y}
& M & 1&&6 \\
 & F & 2&&4 \\ \hline
\multirow{2}{*}{N}
& M & 4&&2 \\
 & F & 6&&1 \\
\hline
\end{tabular}
\end{center}

We observe \begin{align} \label{neqn1} P(Y|M,H) &=
\frac{1}{5}<\frac{2}{8}=P(Y|F,H);\\
P(Y|M,G) &= \frac{6}{8}<\frac{4}{5}=P(Y|F,G).
\label{neqn2}
\end{align}
That is, departmentwise women applicants are favored and hence there is no bias against them.\\
But, considering   the marginal table for $A$ and $X$,

\begin{center}
\begin{tabular}{|cc|lcr|}
\hline &&&X&\\
& &M & &F \\
\hline
 \multirow{2}{*}{A}
& Y & 7&&6 \\
 & N & 6&&7 \\
\hline
\end{tabular}
\end{center}

we get
\begin{equation}\label{eqn3sp} P(Y|M) =
\frac{7}{13}>\frac{6}{13}=P(Y|F),
\end{equation}
showing that overall men do better than women.
 The reversal of the
inequality in \eqref{eqn3sp}, contrary to the ones in  \eqref{neqn1}-\eqref{neqn2},   is called Simpson's paradox.\\

 Why does this paradox occur? To answer this question, we need to look for additional information the data contains.
First, look at the marginal table between $A$ and $D$:

\begin{center}
\begin{tabular}{|cc|lcr|}
\hline &&&D&\\

& &H & &G \\
\hline
 \multirow{2}{*}{A}
& Y & 3&&10 \\
 & N & 10&&3 \\
\hline
\end{tabular}
\end{center}
We have
\begin{equation*}\label{eq5} P(Y|H) =
\frac{3}{13}<\frac{10}{13}=P(Y|G),
\end{equation*}
showing that getting admission in history is tougher than in geography.

 Next,  look at the marginal table between $X$ and $D$:

\begin{center}
\begin{tabular}{|cc|lcr|}
\hline &&&D&\\
& &H & &G \\
\hline
 \multirow{2}{*}{X}
& M & 5&&8 \\
 & F & 8&&5 \\
\hline
\end{tabular}
\end{center}
It is clear that
\begin{equation*}\label{eq4} P(H|M) =
\frac{5}{13}<\frac{8}{13}=P(H|F),
\end{equation*}
because more women applied to history department than in
geography.
Since the admission in history is tougher than in geography, the reversal in
\eqref{eqn3sp} occurs. }
\end{example}

 A probabilistic justification due to Blyth (1972) is the following.
Note that
\begin{eqnarray*}
P(Y|M)&=& P(Y|M,H)P(H|M)+P(Y|M,G)P(G|M) =
\frac{1}{5}\cdot\frac{5}{13}+\frac{6}{8}\cdot\frac{8}{13}=\frac{7}{13}\\
     &=&E_{D|M}P(Y|M,D)\\
P(Y|F)&=& P(Y|F,H)P(H|F)+P(Y|F,G)P(G|F) = \frac{2}{8}\cdot
\frac{8}{13}+\frac{4}{5}\cdot \frac{5}{13}=\frac{6}{13}\\
    &=& E_{D|F}P(Y|F,D)
\end{eqnarray*}
 Because of $\frac{5}{13}=P(H|M) < P(H|F) = \frac{8}{13} $,
the reversal  $P(Y|M)>P(Y|F)$ occurs .
Observe  that the conditional distribution $\mathcal{L}(D|X)$ also plays a key role in
Simpson's paradox.

\begin{example} {\em
Consider the following data (Agresti, (1990)) concerning death
penalty (D),  race of the accused (A) and  race of the victim (V).
Also, let $W$ and $B$ denote the white and the black,
respectively.

\begin{center}
\begin{tabular}{|c|c|lcr|}
\hline
& & & D& \\
\hline
 A&V&Y&&N\\
 \hline
\multirow{2}{*}{W}
& W & 19&&132 \\
 & B & 0&&9 \\ \hline
\multirow{2}{*}{B}
& W & 11&&52 \\
 & B & 6&&97 \\
\hline
\end{tabular}
\end{center}

 First look at the marginal table corresponding to $A$ and $D$. Let $A_{W}$ or $A_{B}$ denote that  accused is a  $W$ or $B$
(similar meaning for $V_{W}$ and $V_{B}$).

\begin{center}
\begin{tabular}{|cc|lcr|}
\hline &&&D&\\
& &Y & & N \\
\hline
 \multirow{2}{*}{A}
& W & 19&& 141 \\
 & B & 17 &&149 \\
\hline
\end{tabular}
\end{center}

\noindent From the above table,
\begin{center}
$P(Y|A_{W})  =  19/160  =   0.12 >  17/166  = 0.10 =P(Y|A_{B})$,
\end{center}
 showing that the white accused are more likely to get death penalty.
However, if we consider the victim's race also, we have from the first table that
the association is reversed for both black and white victims. For,
\begin{center}
$P(Y|A_{W}V_{W})$  =  19/151  =   0.126  $<$ $P(Y|A_{B}V_{W})$ =
11/63  =   0.175
\end{center}
Also,
\begin{center}
$P(Y|A_{W}V_{B})  =  0 <  P(Y|A_{B}V_{B})  = 6/103 = 0.058.$
\end{center}
Thus, Simpson's paradox occurs.}
\end{example}
\subsection{Simpson's paradox for  events}
Blyth (1972) first gave the probabilistic interpretation of
Simpson's paradox, in terms of conditional probabilities. It may happen that for three
events A, B and C,
\begin{equation}\label{eqn66n}
P(A|B) < P(A|B^{c}),
\end{equation}
while
\begin{equation}\label{eqn77n}
P(A|BC) >  P(A|B^{c}C), \quad
P(A|BC^{c}) >  P(A|B^{c}C^{c}).
\end{equation}
As the inequalities in (\ref{eqn77n}) are reversed, compared to
(\ref{eqn66n}), the Simpson's paradox occurs.

Since $P(A|B)$ and $P(A|B^{c})$ are  the following weighted
averages, namely,
\begin{eqnarray}\label{eqn8n}
P(A|B) & = & P(A|BC)P(C|B) + P(A|BC^{c})P(C^{c}|B), \\
P(A|B^{c}) & = & P(A|B^{c}C)P(C|B^{c}) +
P(A|B^{c}C^{c})P(C^{c}|B^{c}), \nonumber
\end{eqnarray}
he pointed out that the reversal happens because the weights $P(C|B)$
and $P(C^{c}|B)$ for  $P(A|B)$ are different than the weights
$P(C|B^{c})$ and $P(C^{c}|B^{c})$ for $P(A|B^{c})$.
 Note if  $B$ and $C$ are independent, then the weights for
$P(A|B)$ and $P(A|B^{c})$ are equal and  the inequalities of (\ref{eqn77n}) will carry over to
 $P(A|B)$ and $P(A|B^{c})$ also.  In other words, the Simpson's
paradox can not happen  in this case.
Thus, Simpson's paradox occurs because of the association between $B$
and $C$.

\begin{remark} {\em Look at Example 2 again. The marginal
table for the  race of the victims (V) and the  race of the accused (A) is:

\begin{center}
\begin{tabular}{|cc|lcr|}
\hline &&&V&\\
& &W & & B \\
\hline
 \multirow{2}{*}{A}
& W & 151 && 9 \\
 & B & 63 && 103 \\
\hline
\end{tabular}
\end{center}
 From the above table, the conditional probabilities for the events $V$ and $A$ are:
\begin{eqnarray}
P(V_{W}|A_{W}) & = & 0.94 \nonumber; ~~P(V_{W}|A_{B})  =  0.38 \nonumber \\
P(V_{B}|A_{W}) & = & 0.06 \nonumber;~~ P(V_{B}|A_{B})  =  0.62
\nonumber
\end{eqnarray}
showing that there is a strong (marginal) association between $V$ and $A$ and
 leading to Simpson's
paradox.}
\end{remark}

\begin{remark} {\em
It is well known that the genesis of Simpson's paradox lies in a simple fact about
proportions.
 There exist positive integers such that
$\displaystyle\frac{k}{l} < \displaystyle\frac{K}{L}$ and
$\displaystyle\frac{m}{n} < \displaystyle\frac{M}{N}$,
but  $\displaystyle\frac{k + m}{l + n} >
\displaystyle\frac{K + M}{L + N}$.
For example, $\displaystyle\frac{1}{6} <
\displaystyle\frac{2}{9}$  and $\displaystyle\frac{5}{7} <
\displaystyle\frac{3}{4}$,
but
$\displaystyle\frac{6}{13} > \displaystyle\frac{5}{13}$ .
This explains why Simpson's paradox occurs in the analysis of some
contingency tables.}
\end{remark}

\subsection{Marginal versus conditional association}
   For the analysis of contingency tables,  Lindley and Novick (1981)
  argue that there is no statistical
criterion that would guard against drawing wrong conclusions or
would indicate which table (conditional or marginal) represents
the correct answer. However, they suggested that if $C$ is
influenced by $B$,  then $C$ should not be treated as a
confounding variable. Pearl (1995)  also suggested that if $C$ is
affected by $B$, then marginal table, rather than the conditional
ones, should be used for inference. Thus, causal considerations
must be used along with inference.
However, there are other researchers who argue Simpson's paradox should not be
viewed in terms of causality, as the reversal is real and is not
causal. Hence, the paradox is a statistical phenomenon that can be
analyzed and  avoided using
tools of statistical techniques.

 One way to avoid Simpson's would be to use a
 randomized experiment which is not always feasible. Cornfield
{\it et al.}~(1959) proposed the minimum effect size criterion to explain  an observed
association measure $\rho (A, B)$ between $A$ and $B$,  if it is
spurious.
 If $B$ has no effect  or less effect than that of $C$   on the
likelihood of $A$, then we would expect
\begin{equation*}\label{eqn4}
\frac{P(C|B)}{P(C|B^{c})} > \frac{P(A|B)}{P(A|B^{c})},
\end{equation*}
or the risk difference condition, namely,
\begin{equation*}\label{eqn5}
P(A|C) - P(A|C^{c}) \geq P(A|B) - P(A|B^{c}).
\end{equation*}
 Schield (1999) suggested that this condition could be used as a
simple method for deciding whether $C$ has the strength-the effect
size necessary-to reverse the association $ \rho(A, B)$.

\subsection{ Simpson's paradox as an
association reversal phenomena }

  Samuels (1992) showed that Simpson's paradox  between events can be viewed as a particular case of
  association reversal   phenomena for random
variables/distributions, which we describe now.
Let $(Y, X, W) \sim F$  and for example $F_{X}(x)$ denote the marginal distribution of $X$.  We say $W$ is not doubly linked to
$(Y, X)$ if at least one of the following condition holds:
$$(a)~ W \perp Y,  ~ (b)~ W\perp X,
~(c)~  W \perp Y|X,  ~(d)~  W \perp X|Y. $$
Otherwise, it is doubly linked to $(Y, X)$, {\it i.e.}, $W$ is linked to
both $Y$ and $X$. Henceforth, $ X \perp Y|W $ denotes the conditional independence of $X$ and $Y$, given $W$.

 An  association reversal can be defined for  any relation $R= R(Y, X)$ which
denotes the directional association between any two random variables $X$ and $Y$. Henceforth,
$\uparrow$ and $\downarrow$ means respectively nondecreasing and nonincreasing. Some  relations  studied in
the literature are:
\begin{namelist}{xxx}
\item ${\mathcal{R}_{1}}$: (stochastically increasing ) $Y \uparrow X$ if $P(Y
> y | X = x)$  is $\uparrow$ in $x$ for all $y$
.

\item ${\mathcal{R}_{2}}$: (mean incresing) $Y \uparrow X$ if $E(Y | X = x)$ is
$\uparrow$ in $x$  .

\item ${\mathcal{R}_{3}}$: (positive quadratic
dependence) $Y \uparrow X$ if $F (y, x)  \geq
F_Y(y)F_X(x)$ for all ($x, y$)

\item ${\mathcal{R}_{4}}$: (covariance increasing) $Y \uparrow X$ if $Cov (X, Y) >$ 0 .
\end{namelist}
The relations ${\mathcal{R}_{3}}$ and ${\mathcal{R}_{4}}$ are
symmetric in $Y$ and $X$. The above relations can also be defined
for $\downarrow$ case also.

  Samuels (1992) proved that the joint distribution  $F$ cannot exhibit association reversal with respect to $ {\mathcal{R}_{3}}$ if $W$ is not doubly
linked to $(Y, X)$, that is, when one of the conditions (a) to (d) is true. But,
 the above result is not true for the relation ${\mathcal{R}_{4}}$. For
example, the condition $W \perp X|Y$ is not sufficient to prevent the association reversal for $
{\mathcal{R}_{4}}$.  However, $W \perp Y$ prevents association reversal for several $\mathcal{R}$'s. See Samuels (1992) for some additional results in this direction.

\subsection{Linear regression models }
 Let
\begin{equation}\label{eqn6}
E(Y|X, W)  =  \beta_{0} + \beta_{1}X + \beta_{2}W
\end{equation}
with $\beta_{1} \leq$ 0. Also, let $\eta = \beta_2 Cov(X, W).$
Samuels (1992) showed that the distribution $F$ exhibits positive association reversal for
${\mathcal{R}_{4}}$ iff $ \eta
>$ 0, and $|\eta| > |\beta_{1}| Var(Y).$
As  a corollary, the association reversal with respect to ${\mathcal{R}_{2}}$ holds.

 Suppose the  marginal model is also linear defined by
\begin{equation*}\label{eqn7}
E(Y|X)  =  \tilde{\beta_{0}}  +  \tilde{\beta_{1}}X,
\end{equation*}
where
\begin{equation*}
\tilde{\beta_{1}}  =  \frac{Cov(Y, X)}{Var(X)} .
\end{equation*}
Note from (\ref{eqn6}),
\begin{equation*}
\beta_{1}  =  \frac{Cov(Y, X|W)}{Var(X|W)}.
\end{equation*}

 It is possible that $\beta_{1} <$ 0 while $\tilde{\beta_{1}}
>$ 0, implying the occurrence of Simpson's paradox for the regression coefficients.
Some sufficient conditions for the Simpson's paradox have been recently
discussed in Chen, Bengtsson and Ho (2009).

\subsection {Simpson's paradox in survival analysis}

  In the context of
the survival analysis, it is possible that increasing the value of
a covariate $X$ has a positive effect on a failure time $T$, but
this effect may be reversed  when conditioning on another possible
covariate $Y$. When studying causal effects and influence of
covariates on a failure time $T$, this aspect appears paradoxical
and creates suspicion on the real effect of $X$. These situations
 may be seen as a kind
of Simpson's paradox.

\vone Let  $X$ and $Y$ be the covariates having effect on $T$.
Simpson's paradox occurs in survival probability at $(t, s)$ if
\begin{eqnarray}\label{eqn8}
\begin{array}{lll}
(a) & P(T > t + s | T > t, X = x, Y = y) \downarrow  \textrm{ in
$x$, for all $y$ and}  \\
(b) & P(T > t + s | T > t, X = x) \uparrow  \textrm{ in $x$}.
\end{array}
\end{eqnarray}
 Scarsini and Spizzichino (1999) discussed the Simpson's
paradox for different notions of positive dependence and aging.

 Let $h(t|x)$ denote the conditional hazard rate, given $X=x$, defined by
\begin{equation*}
h(t|x)  = \lim_{s \downarrow 0} P(t < T < t + s | T > t; X = x).
\end{equation*}
Then, Simpson's paradox  for the hazard rate occurs if
\begin{equation} \label{eqn9}
h(t|x, y) \uparrow  \textrm{in}~ x,~ \textrm{for all}~ y,~ \textrm{but} ~
h(t|x) \downarrow  \textrm{in}~ x.
\end{equation}
 Di Serio, Rinott and Scarsini (2009) showed that Simpson's paradox
occurs naturally in the context of survival analysis. They studied
the range ($t, s$) for which (\ref{eqn8}) holds, and showed that
under certain
conditions it holds for all $t, s >$ 0.
They discussed Simpson's paradox  for the linear
transformation model defined by
\begin{equation}\label{eqn10}
K(T)  =  - \beta_{x}X  -  \beta_{y}Y  +  W, \end{equation} where
$K$ is increasing and $W \perp (X, Y)$.

 Suppose the covariates $X$ and $Y$ satisfy the model
\begin{equation*}
Y  =  \eta(X)  +  V,
\end{equation*}
where $\eta$ is $\uparrow$ and $X \perp V$.
 That is, $(Y|X = x)$ has density
\begin{equation}\label{eqn11}
f(y|x)  =  f_{v} (y  - \eta(x)).
 \end{equation}

\noindent Their main result is the following:

\begin{theorem}
Let $T$ follow the  model (\ref{eqn10}) and the conditional distribution $\mathcal{L}(Y|X)$
follow
(\ref{eqn11}). Assume $V$ and $W$ have increasing failure rates ($IFRs$). Then, \\
(i) Simpson's paradox in survival probability  defined in (\ref{eqn8})
occurs for all $t, s  >$  0 if
\begin{equation}\label{eqn12}
\beta_{y} < 0 < \beta_{x} \hspace{4mm} and \hspace{4mm} \beta_{x}x
+ \beta_{y}\eta(x)~ \mbox{is decreasing in} ~ x.
\end{equation}
(ii) If $W$ and $V$ are both strictly $IFRs$, then Simpson's
paradox  for hazard rate in (\ref{eqn9}) occurs for all $t, s  >$
0 if and only if  (\ref{eqn12}) holds.
\end{theorem}
For example, when $(Y|X = x) \sim N(\mu + \rho x, 1 - \rho^{2})$,
then (\ref{eqn9}) holds
$\Longleftrightarrow \beta_{y} < 0 < \beta_{x}$ and $\rho > \displaystyle\frac{\beta_{x}}{|\beta_{y}|}.$

  Note that the linear transformation model corresponds to the
survival function
\begin{equation*}
\overline{F}_{T} (t|x, y)  =  \overline{F}_{W} ( K(t)  +  \beta_{x}x  + \beta_{y}y).
\end{equation*}
 It can be seen that the Cox's  (1972) {\it proportional hazard model}
\begin{equation*}
h(t|x, y)  =  h_{0}(t)e^{( \beta_{x}x + \beta_{y}y)},
\end{equation*}
where $h_{0}(t) >$ 0, is a special case of linear transformation
model with
$\overline{F}_{W}(t)  =  -e^{t}, ~ t \hspace{1mm} \epsilon \hspace{1mm} \mathbb{R}$.
Also, the {\it proportional odds model} (Pettitt (1984))
corresponds to the case
$\overline{F}_{W}(t)  =  1/({1  +  e^{t}}), ~ t \hspace{1mm} \epsilon \hspace{1mm} \mathbb{R},$
the {\it logistic}  distribution.


 It is interesting to note that even when covariates $X$ and $Y$ are independent, there exists the
choice of parameters in Cox model for which Simpson's paradox in
survival probability occurs for some $t, s >$ 0. See Di Serio {\it et al.} (2009) for more details.
However, the association or the dependence between $X$
and $Y$, modeled  through conditional distributions, is the main source of Simpson's paradox.

\subsection {Simpson's paradox for an association/dependence measure}
 Kendall's $\tau$ and Spearman's $\rho$ are well known
measures of concordance, a certain form of dependence.
For example, the dependence of $Y$ on $X$ is called
stochastically increasing if $P(Y >y\mid X= x)$ is increasing in
$x$ for all $y.$
 In other words,  when $X$ is continuous and  the partial
derivative exists (Cox and Wermuth (2003)), the conditional distribution function $F(y|x)$ satisfies
\begin{eqnarray}
\frac{\partial{F(y\mid x)}}{\partial x} \leq 0, \label{eqn13}
\end{eqnarray}
for all $y$ and $x$, with strict inequality in a region of
positive probability.

 Let $W$ be a covariate. Then,
\begin{eqnarray} \label{eqn14dd}
\frac{\partial{F(y\mid x)}}{\partial x} & = & \displaystyle\int
\frac{\partial{F(y\mid x, w)}}{\partial x}
 f(w\mid x) \,dw
  + \int F(y\mid x,w) \frac{\partial{f(w\mid x)}}{\partial x} \,dw.
\end{eqnarray}

\vone \noindent If $W\perp X$, then
$\displaystyle \frac{\partial{f(w\mid x)}}{\partial x}= 0$
and hence (Cox (2003))
\begin{equation*}
\label{eqn15} \displaystyle\frac{\partial{F(y\mid x)}}{\partial
x}= \displaystyle\int \frac{\partial{F(y\mid x, w)}} {\partial x}
f(w)dw,
\end{equation*}
showing that
 $$\displaystyle\frac{\partial{F(y\mid x, w)}}{\partial x}\leq 0
\Longrightarrow \displaystyle\frac{\partial{F(y\mid x)}}{\partial
x}\leq 0, \ {\rm{ for ~ all}}\
 y,x \ {\rm and} \ w. $$
Thus, $Y$ remains stochastically increasing in $x$ after
marginalization over the covariate $W.$
 Note in general
   (see (\ref{eqn14dd})) it is possible that
$$\displaystyle \frac{\partial{F(y\mid x, w)}}{\partial x}\leq 0,
\mbox{~for all~} (y, x, w), ~\mbox{but}~
 \displaystyle \frac{\partial{F(y\mid x)}}{\partial x}>0$$ for some $y$ and $x$, implying
 Simpson's paradox.

\section{Collapsibility}

Collapsibility is a concept closely related to that of Simpson's paradox. Generally, whenever Simpson's
paradox does not occur, the collapsibility issue arises naturally as a dimension
reduction problem in the context of data analysis. It was originally associated with the analysis of  contingency
tables and so we start with the same.

\subsection{Collapsibility of contingency tables}

  Let $ X_1,\cdots,X_n $ be a set of n categorical variables, where
$ X_j \in \{1,\cdots,m_j\}, ~ 1 \leq j \leq n $.
Let $  i=(i_1,\cdots,i_n) $ denote a cell of the $n$-dimensional table, and $p(i)$ denote the cell probability
with $p(i)>0$ and $\displaystyle\sum_i p(i)=1$. \\
Let $\bar n=\{1,2,\cdots,n\}$.
Define $l^{(n)}(i)=l^{(n)}(i_1,\cdots,i_n)=\ln p(i_1,\cdots,i_n)$.
Let $A=(a_1,\cdots,a_r), a_j\in \bar n$, $i_A=(i_{a_1},\cdots,i_{a_r})$, and $|A|$
denotes the cardinality of A.
Define, for any subset $A \subset \bar{n}$,
\begin{eqnarray*}
l_{A}^{(n)}(i_A)=\displaystyle\sum_{i_j:j\in A^c} l(i_1,\cdots,i_n); ~~\tilde{l}_{A}^{(n)}(i_A)=\displaystyle{\frac{1}{\displaystyle\prod_{j\in A^c}m_j}}l_A(i_A).
\end{eqnarray*}

 For the $n$-dimensional table, let
\begin{eqnarray}
l^{(n)}(i)= \displaystyle\sum_{Z\subseteq \bar{n}} \tau^{(n)}_Z(i_Z)
\label{loglin}
\end{eqnarray}
be the log-linear model (LLM),  where $\tau^{(n)}_Z(i_Z)$ is the $r$-factor
 interaction parameter when $|Z|=r$. Then, it can be seen that (Vellaisamy and Vijay (2007))
\begin{eqnarray}
\tilde{l}_A^{(n)}(i_A)=\displaystyle\sum_{Z\subseteq A}\tau_Z^{(n)}(i_Z), ~~~~\forall A\subseteq \bar{n}.
\label{logeq}
\end{eqnarray}

 For example, when n=3,
\begin{eqnarray*}
l^{(3)}(i_1,i_2,i_3)&=& \tau_{123}^{(3)}(i_1,i_2,i_3)+ \tau_{12}^{(3)}(i_1,i_2)+\tau_{13}^{(3)}(i_1,i_3)+
\tau_{23}^{(3)}(i_2,i_3)\nonumber\\
&& + \tau_{1}^{(3)}(i_1)+\tau_{2}^{(3)}(i_2)+
\tau_{3}^{(3)}(i_3)+\tau_{\phi}^{(3)}\\
&=& \displaystyle\sum_{A}\tau_A^{(3)}(i_A),\nonumber
\end{eqnarray*}
where $A$ is any subset of $\{1,2,3\}$.
Then
\begin{eqnarray*}
\tilde{l}_{12}^{(3)}(i_1,i_2) &:=&\frac{1}{m_3}\displaystyle\sum_{i_3}l^{(3)}(i_1,i_2,i_3), \nonumber\\
                     &=&\tau_{12}^{(3)}(i_1,i_2)+ \tau_{1}^{(3)}(i_1)+\tau_{2}^{(3)}(i_2)+\tau_{\phi}^{(3)}\nonumber\\
                     &=&\displaystyle\sum_{Z\subseteq \{1,2\}} \tau_Z^{(3)}(i_Z).
\end{eqnarray*}
Indeed, the interaction factor admits the following representation
\begin{eqnarray}
\tau_A^{(n)}(i_A)
=\sum_{Z\subseteq A}(-1)^{|A-Z|}\tilde{l}_Z^{(n)}(i_Z), ~~~\forall A\subseteq \bar{n}.
\label{oeq3}
\end{eqnarray}
 Whittemore (1978) defined first $\tau_A^{(n)}(i_A)$ as a straightforward extension
and  remarked later
that $l(i)=\displaystyle\sum_{Z\subseteq \bar{n}}\tau_Z(i_Z)$.
 Recently, Vellaisamy and Vijay (2007) gave a direct proof of \eqref{oeq3} by considering $\tilde{l}_A$ as the function on the poset $(\cal P, \subseteq)$,
and using M\"obius inversion theorem. Note also that, by M\"obius inversion theorem,
(\ref{logeq}) holds iff (\ref{oeq3}) holds.

 Let now for simplicity  $A=\{1,\cdots, r\}$, and $B=\{1,\cdots,r, {r+1},\cdots, s\}$, where $r\leq s < n$.
Define $p_B(i_B)= \displaystyle\sum_{i_j:j\in B^c}p(i)$,
the cell probabilities of the marginal (condensed over $B^c$) table. Define, similarly,
\begin{eqnarray}
l^{(s)}(i)= \ln (p_B(i_B)) = \displaystyle\sum_{Z\subseteq B}\eta^{(s)}_Z(i_Z)
\end{eqnarray}
be the LLM for the marginal table. Then as seen in the LLM for the full table,

\begin{eqnarray}
\tilde{l}_Z^{(s)}(i_Z)=\frac{1}{\displaystyle\prod_{j\in B\setminus Z}m_j}
\displaystyle\sum_{i_j:j\in B\setminus Z}l^{(s)}(i) = \displaystyle\sum_{A\subseteq Z}\eta^{(s)}_A(i_A),
\label{defls}
\end{eqnarray}
for any  $Z\subseteq B$. The following definition of collapsibility is due to  Whittemore (1978).
\begin{definition}  An n-dimensional table is said to be  collapsible
into an s-dimensional table  with
respect to $\tau_A^{(n)}$, $A\subseteq B$, if
$ (i)~ \tau_A^{(n)}=\eta_A^{(s)}$ and strictly collapsible if, in addition to $(i)$,
$(ii) \tau_Z^{(n)}=0, ~~\forall ~Z\supseteq A, Z\cap B^c\neq \phi$ holds.
\end{definition}

Let
\begin{eqnarray}
d^{(B)}(i_B)=l^{(s)}(i)- \tilde{l}_B^{(n)}(i_B)
\label{defd}
\end{eqnarray}
and for any $Z\subseteq B$
\begin{eqnarray}
\tilde{d}^{(B)}_Z(i_Z)= \frac{1}{\displaystyle\prod_{j\in B\setminus Z}m_j}
\sum_{i_j:j\in B\setminus Z}d^{(B)}(i_B).
\label{dz}
\end{eqnarray}

The next result (Vellaisamy and Vijay (2007)) characterizes the conditions for collapsibility.
\begin{theorem} Let
$\delta_Z=(\eta^{(s)}_Z-\tau^{(n)}_Z)$, for $Z\subseteq B$. An $n$-dimensional table is collapsible
to an $s$-dimensional table with respect to $\tau_{A}^{n}$ if and only if
\begin{eqnarray}
 \tilde{d}^{(B)}_A(i_A)=\displaystyle\sum_{Z\subset A}\delta_Z(i_Z) \Longleftrightarrow
\displaystyle\sum_{Z\subseteq A}(-1)^{|A-Z|} \tilde{d}^{(B)}_Z(i_Z)=0,
\end{eqnarray}
where $\tilde{d}_Z^{(B)}$ is defined in (\ref{dz}), and $A\subseteq B$.
\label{th1}
\end{theorem}

 We next mention an important result for the strict collapsibility for hierarchical log-linear models (HLLM), a subclass of LLMs, defined as follows:

\begin{definition}
A LLM $l^{(n)}(i)=\displaystyle\sum_{Z\subseteq \bar{n}}\tau_Z^{(n)}$
is said to be hierarchical if $\tau_B^{(n)}\neq 0~\Longrightarrow~~ \tau_A^{(n)}\neq 0$ for $A\subset B$
or equivalently $\tau_C^{(n)}=0 \Longrightarrow~~ \tau_D^{(n)}= 0$ for $D\supset C$.
\end{definition}

 Let now $ \bar{n}= A + B + C$.  For a HLLM, Bishop, Fienberg and Holland (1975) (BFH (1975)) showed that the n-dimensional table is collapsible
into a s-dimensional table (over C) with respect to $\tau^{(n)}_{A\cup V}$,
where $V\subseteq B$, iff $\tau_Z^{(n)}=0$, for all $Z\cap A\neq \phi$ and $Z\cap B\neq \phi$,
that is, $X_A\perp X_C|X_B$.
 Later, Whittemore (1978) showed that they are
only sufficient but not necessary.
Recently, Vellaisamy and Vijay (2007) showed that those conditions are necessary and sufficient for
strict collapsibility with respect to a set of interaction parameters, which is stated below.
\begin{theorem} Let $\bar{n}= A + B + C$ be such that $|A\cup B|=s$ and
$|C|=n-s$. Then, an n-dimensional table is strictly collapsible (over C) into an s-dimensional
table with respect to the set $C_L=\{\tau_L| L\subseteq A\cup B; L\cap A\neq \phi\}$
if and only if $X_A\perp X_C|X_B$.
\end{theorem}

 It follows from the above result that for $k\in \{1,2\}$,  a 3-dimensional table is strictly collapsible into a 2-dimensional table with
respect to $\tau_k^{(3)}$ and $\tau_{12}^{(3)}$ iff
$\tau_{123}^{(3)}=0$ and $\tau_{k3}^{(3)}=0$.
Note also  that when k=1, the conditions $\tau^{(3)}_{123}=0$ and $\tau^{(3)}_{13}=0$
are nothing but BHF's (1975) sufficient conditions for collapsibility with
respect to $\tau_{12}^{(3)}$ or $\tau_{23}^{(3)}$.

For some recent  results on the collapsibility of full  tables based on conditional tables/models, one may refer to
 Vellaisamy and Vijay (2010). The concept of collapsibility for contingency tables was later extended to the study of regression models by Wermuth (1989) and
several others.

\subsection{Collapsibility of regression coefficients}

Let Y be a continuous response variable, X be a continuous influence variable and A be a discrete
(background) variable with levels $i=1,2,\cdots,I$. Initially, the problems of collapsibility were
 addressed only for  parallel regression models
 (Wermuth (1989)) defined by
\begin{eqnarray}
E(Y|X=x,A=i)= \alpha_{yx}(i) + \beta_{yx} x,~~~~1\leq i\leq I,
\label{1}
\end{eqnarray}
where $\beta_{yx}= \displaystyle\frac{\sigma_{yx}(i)}{\sigma_{xx}(i)}
=\displaystyle \frac{Cov(Y,X|A=i)}{V(X|A=i)}$
is the regression
coefficient. Since $\sigma_{yx}(i)=\beta_{yx}\sigma_{xx}(i)$ for all $i$, we have
\begin{eqnarray*}
\sigma_{yx}(A)=Cov(Y,X|A)=\beta_{yx}V(X|A).
\label{b2}
\end{eqnarray*}

Let us now introduce the following notation:
\begin{eqnarray*}
&&\mu_y(A)=E(Y|A),~~~~ \mu_x(A)=E(X|A),~~~~\sigma_{xx}(A)=V(X|A),\\
&&\sigma_{yx}=Cov(Y,X),~~~~\sigma_{xx}= V(X), P(A=i)=\pi_i>0. \nonumber
\label{a}
\end{eqnarray*}
In general, $\beta_{yx}(A)=\displaystyle ({\sigma_{yx}(A)}/{\sigma_{xx}(A)})$
is a function of A. In the parallel regression model, $\beta_{yx}(A)=\beta_{yx}$. Note also from the
model \eqref{1},
$\mu_y(A) = E_{X|A}(\alpha_{yx}(A) + \beta_{yx} X).$
 The following definition of collapsibility is due to Wermuth (1989).

\begin{definition} The parallel regression coefficient $\beta_{yx}$ is said to be collapsible
over $A$ if $\beta_{yx}=\tilde{\beta}_{yx}$, where $\tilde{\beta}_{yx}$ is the regression coefficient
for the marginal linear model
\begin{eqnarray}
E(Y|X=x)= \tilde{\alpha}_{yx} + \tilde{\beta}_{yx} x.
\label{d1}
\end{eqnarray}
\end{definition}

 The above model implies $\sigma_{yx}=\tilde{\beta}_{yx}\sigma_{xx}$.
 Next we present a necessary and sufficient condition for collapsibility, due to  Wermuth (1989) (see also Vellaisamy and Vijay (2008) for a simple probabilistic proof).

\begin{theorem} \label{t1}  The regression coefficient
$\beta_{yx}$  of the parallel regression model {\rm(\ref{1})} is collapsible over $A$ if and only if
$Cov_A(\alpha_{yx}(A), \mu_x(A))=0.$
\end{theorem}

 As a corollary, $(i)~ \alpha_{yx}(A)$ or $\mu_x(A)$ is degenerate, or
 $ (ii)~ \beta_{yx}= ({\mu_y(A)}/{\mu_x(A)}) $ {\it a.e.}, with $\mu_x(A)\neq 0$, is a sufficient condition for collapsibilty.

The next  result, due to Vellaisamy and Vijay (2008),  is more general, as it does not assume $(Y,X,A)$
 follows a conditional Gaussian distribution, a condition usually assumed in the literature.

\begin{theorem} \label{nn} The regression coefficient $\beta_{yx}$ of the model (\ref{1}) is collapsible if
\begin{itemize}
\item[{\rm(i)}] $Y\perp A|X$ ~~or
\item[{\rm(ii)}] $X\perp A|Y$ ~ and ~ $\displaystyle {V_A(\mu_y(A))} {E_A(\sigma_{yy}(A))}
=\displaystyle {V_A(\mu_x(A))} {E_A(\sigma_{xx}(A))}.$
\end{itemize}
\end{theorem}

\begin{remark} Note that if $X\perp A|Y$ and $Y\perp A$, then $A\perp(X,Y)$ (Whittaker (1990)) and hence
condition (ii) of the above result is satisfied. Thus, $\beta_{yx}$ is collapsible.
\end{remark}

\subsection{Random Coefficient Regression Models}
 The  condition that ${Cov((Y,X)|A=i)}/{V(X|A=i)}$ is independent of the levels of A
is stronger and  may not hold in several real-life applications. For example,
 the well-known  degradation models of the form
\begin{eqnarray}
\label{deg}
y_i(t)= \alpha(i)-\beta(i)t,
\end{eqnarray}
where $y_i(t)$ denotes the log-performance of specimen $A=i$ as a function of
age t, shows that different specimens have different linear degradation.
Such models arise in accelerated life-testing problems (Nelson (2004), p. 530).
The model (\ref{deg}) can be written in  the form
\begin{eqnarray*}
E(Y|X=x,A=i) = \alpha_{yx}(i) + \beta_{yx}(i)x,
\end{eqnarray*}
or equivalently
\begin{eqnarray}
E(Y|X, A) = \alpha_{yx}(A) + \beta_{yx}(A)X,
\label{ac11}
\end{eqnarray}
where  $\beta_{yx}(A)= (\sigma_{yx}(A)/  \sigma_{xx}(A)) $, is a random coefficient regression model. For these models,
Vellaisamy and Vijay (2008) introduced and studied average collapsibility which we discuss next.
\begin{definition} The random regression coefficient $\beta_{yx}(A)$
is said to be  average collapsible (A-collapsible) if
$\tilde{\beta}_{yx}= E_A(\beta_{yx}(A))$,
where $\tilde{\beta}_{yx}$ is the regression coefficient of the marginal linear model
\begin{eqnarray}
\label{3}
E(Y|X=x) = \tilde{\alpha}_{yx} + \tilde{\beta}_{yx}x.
\end{eqnarray}
\end{definition}

 The following result  generalizes Theorem \ref{t1}.

\begin{theorem} \label{avt}
The random regression coefficient $\beta_{yx}(A)$ of the model {\rm(\ref{ac11})} is A-collapsible if and only if
\begin{eqnarray}
E_A(\beta_{yx}(A))V(\mu_x(A)) = Cov(\beta_{yx}(A), ~\sigma_{xx}(A)) + Cov(\mu_y(A),~ \mu_x(A)).
\label{ac3}
\end{eqnarray}
\end{theorem}

 It is of practical interest to know the conditions under which both the random regression coefficients  are collapsible.
\begin{theorem} \label{tt} Consider the random coefficients regression model {\rm(\ref{ac11})}
with $P(X=0)=0$. Then $\alpha(A)$ and $\beta(A)$ are
both A-collapsible if one of the following conditions holds:
\begin{eqnarray*}
{\text(i)}\quad E(\alpha(A)|X)= \tilde{\alpha}~a.e.,\quad
{\text(ii)}\quad E(Y|X,A)=E(Y|X)~~ a.e.\label{21}
\end{eqnarray*}
\end{theorem}

 Next, we briefly discuss the collapsibility problems for the logistic regression coefficients.
Let Y be a binary response variable taking the values 0 and 1, X be a random vector
of $p$ risk factors and A be a discrete background variable with levels
$i=1,\cdots,I$. Guo and Geng (1995) obtained some collapsibility results for the
regression coefficients of the logistic regression model. We focus only on  random  coefficient
logistic regression model of Y on X, for
the levels of A, defined by
\begin{eqnarray}
\ln \displaystyle\{\frac{P(Y=1|X, A)}{P(Y=0|X, A)}\}= \alpha(A) + \beta^T(A)X.
\label{lg1}
\end{eqnarray}
We say that the logistic regression coefficient vectors $\beta(A)$ is A-collapsible
if $E_A(\beta(A))=\tilde{\beta}$, where $\tilde{\beta}$ is the regression coefficient vector for
the marginal regression model
\begin{eqnarray*}
\ln \displaystyle\{\frac{P(Y=1|X)}{P(Y=0|X)}\}= \tilde{\alpha} + \tilde{\beta}^TX.
\label{lg2}
\end{eqnarray*}
\begin {theorem} \label{tt3}
Let X be a continuous random vector. Then, for the model {\rm(\ref{lg1})},
\begin{itemize}
\item[{\rm(i)}] $A\perp Y|X$ implies $\alpha(A)$ and $\beta(A)$ both are A-collapsible.
\item[{\rm(ii)}] $A\perp X|Y$ implies $\beta(A)$ is A-collapsible.
\end{itemize}
\end{theorem}

 For a proof, see Vellaisamy and Vijay (2008). Finally, we address the collapsibility issues for a certain
measure of dependence for two random variables.

\subsection{Collapsibility of distribution dependence}
Let $F(y|x, w)$ denote the conditional distribution function, where $Y$ is a response variable, $X$ is an
explanatory variable (continuous) and $W$ is a background variable. Then, the
function $ \displaystyle \frac{\partial F(y|x,w)}{\partial x}$, when it exists, is
called a distribution dependence function (Cox and Wermuth (2003)). It  represents
the stochastically increasing property between $X$ and $Y$.  When $X$ is
discrete, the partial differentiation is replaced by differencing
between adjacent levels of $X$.
The following definition is due to Ma, Xie and Geng (2006).

\begin{definition} The distribution dependence function is said to be
homogeneous with respect to $W$ if
$\dfrac{\partial F(y|x,w)}{\partial x}=\dfrac{\partial F(y|x,w')}{\partial x},$
for all $y$, $x$ and $w\neq w'$
and collapsible over $W$ if
$\dfrac{\partial F(y|x, w)}{\partial x}=\dfrac{\partial
F(y|x)}{\partial x},  \ {\rm{ for\  all}}\
 y,x \ {\rm and} \ w $. \nonumber
\end{definition}

 Ma \emph{et al.~}(2006)
showed that the distribution
 dependence function is uniformly collapsible and hence collapsible iff either (a) $Y \perp X | W$; or (b) $X \perp W$ and $\displaystyle
\frac{\partial F(y|x,w)}{\partial x}$ is homogeneous in $w$.
Cox and Wermuth (2003) showed that either condition
(a) or (b) is sufficient to ensure that no  effect reversal  or Simpson's paradox occurs.
 Note that homogeneity is a stronger condition  which may not hold  for
most of the models that are encountered in practice.  For example,  consider a simple
linear regression model defined by
$Y = m(X, W) +\epsilon$,
where $ m(x,w) = \alpha_1x + \alpha_3xw$, and $\epsilon
\sim N(0, 1).$ Let $\phi$ be the standard normal density. Then,
\begin{eqnarray*}
\frac{\partial F(y|x,w)}{\partial x} = -(\alpha_1+ \alpha_3 w)
\phi (y-m(x,w)),
\end{eqnarray*}
and hence is not homogeneous over $W$. For such models, the concept of average collapsibility
introduced  and studied by Vellaisamy (2011) is a very useful concept. Indeed, when the  distribution dependence function
is homogeneous, it reduces to collapsibiltiy.
 We say that the distribution dependence function
$\dfrac{\partial F(y|x, w)}{\partial x}$ is average collapsible over  W
if
\begin{eqnarray*} \label{eq3.1}
E_{W\mid X=x}\left(\dfrac{\partial F(y|x,W)}{\partial x}\right) =
\dfrac{\partial F(y|x)}{\partial x},~\mbox {for all y and x}.
\end{eqnarray*}

 Vellaisamy (2011) showed that  average collapsibility holds if  (i) $Y\perp W\mid X$ or (ii)  $W\perp X$  holds. These conditions are also
necessary when $W$ is a binary variable.
 An example, where average collapsibility holds, follows next.
  Let $\phi(z)$ and $\Phi(z)$ respectively denote the
density and the distribution of $Z \sim N(0, 1).$

\begin{example} \label{ex-lr} {\em
Consider the linear regression model
\begin{equation*}\label{eq2.4n}
 Y= \alpha_1X +\alpha_2 W + \alpha_3 XW + \epsilon,
\end{equation*}
where $\epsilon \perp (X, W)$ and $\epsilon \sim N(0, \sigma^2)$.
Then
    $$ (Y|x, w) \sim N(m(x,w), \sigma^2),$$
where $ m(x, w)=\alpha_1x +\alpha_2 w + \alpha_3 xw$.
Hence,
\begin{equation*}\label{eq2.5n}
 \dfrac{\partial F(y|x, w)}{\partial x}=(\frac{-1}{\sigma})(\alpha_1+ \alpha_3w)
 \phi(\dfrac{y-m(x,w)}{\sigma}),
\end{equation*}
which is not homogeneous.

 Suppose  $W \sim N(0, 1)$  and $W\perp X$.
 Then $(Y|x) \sim N(\alpha_1 x, v^{2}(x, \sigma))$, where $
v^{2}(x, \sigma)= (\alpha_2 +\alpha_3x)^2+ \sigma^2.$
Then, it can be shown that (see Vellaisamy (2011))
$$E_{W\mid X=x}\left(\dfrac{\partial F(y|x,W)}{\partial x}\right)=\dfrac{\partial F(y|x)}{\partial x}, $$
so that average collapsibility  over $W$ holds.}
\end{example}

 Note also from  \eqref{eqn14dd} that average collapsibility
holds if and only if
\begin{eqnarray}\label{eqs1}
\int F(y\mid x,w) \frac{\partial{f(w\mid x)}}{\partial x} \,dw =
0~~ \mbox{for all $(y,x)$}.
\end{eqnarray}

  The conditions (i) and (ii) of average collapsibility are not necessary, unless $W$ is  binary.
  A counter-example follows:

 \begin{example} \label{ex3} {\em
Let $(Y|x, w)\sim U(0,
(x^2+(w-x)^2)^{-1})$ so that
\begin{eqnarray*}
F(y|x,w) = y(x^2+(w-x)^2), ~~0<y<(x^2+(w-x)^2)^{-1}.
\end{eqnarray*}
Assume also $(W|X=x)\sim N(x,1)$ so that
\begin{eqnarray*}
\frac{\partial}{\partial x}f(w|x) = - \phi'(w-x) = (w-x)\phi(w-x).
\end{eqnarray*}
Then it can be seen that (see Vellaisamy (2011))
\begin{equation*}
\int F(y|x,w)\frac{\partial}{\partial x} f(w|x)dw
= 0, ~~\mbox{for all $(y,x)$}.
\end{equation*}
Thus, from (\ref{eqs1}), average collapsibility  over $W$ holds, but
neither condition (i) nor condition (ii) is satisfied.}
\end{example}

\vone
 {\bf Conclusions}.  The examples and the applications discussed in this paper clearly demonstrate that Simpson's paradox is
a crucial aspect in the data analysis and the issue of collapsibility should be looked into only after ascertaining the nonoccurrence of Simpson's paradox. Only recently, the issue of Simpson's paradox for survival
analysis and for certain measures of association has been addressed. However, the conditions of collapsibility for survival models, when the co-variate
is either known or unknown, are yet to be explored. Specifically, the concept of average collapsibility
 (Vellaisamy (2011)) is more
relevant in view of the nature of the Simpson's paradox in survival models (Di Serio {\it et.~al} (2009)). These and other considerations are of practical interest and some of these issues are already under consideration. The findings will be reported elsewhere.

\vone
 {\bf Acknowledgements}. This work was completed while the author was visiting the Department
of Statistics and Probability, Michigan State University, USA. The author is grateful to Professor Hira L. Koul
for  all the  support and encouragement for the timely completion of this work, and for some helpful comments
which improved the presentation of the paper. This research is partially supported by a DST research grant No. SR/S4/MS: 706/10.

\vone
 {\large \bf References}

\begin{namelist}{xxx}

\item Bickel, P. J.,  Hjammel, E. A. and O' Connel, J. W. (1975).
Sex bias in Graduate admissions: Data from Berkeley. {\it
Science}, {\bf 187}, 398-404.

\item Bishop, Y. M. M., Fienberg, S. E. and Holland, P. W. (1975). {\it Discrete Multivariate Analysis:
Theory and Practice}. MIT Press, Cambridge.

 \item Blyth, C.R. (1973). Simpson's paradox and mutually favourable
events. {\it J. Amer. Statist. Assoc.}, {\bf
68}, 746.

\item   Chen, A., Bengtsson, T. and Ho, T.K. (2009). A regression
paradox for linear models: sufficient conditions and relation to
Simpson's paradox. {\it Amer. Statistician}, {\bf 63}, 218-225.

\item   Cornfield, J., Haenszel, W., Hammond, E., Lilienfield, A.,
Shimkin M., {\it et al.} (1959). Smoking and lung cancer: recent
evidence and a discussion of some questions. {\it J.
Nat. Cancer Institute}, {\bf 22}, 173-203.

\item Cox, D. R. (2003). Conditional and marginal association for
binary random variables. {\it Biometrika}, {\bf 90}, 982-984.

 \item  Cox, D.R. and Wermuth, N. (2003). A general condition for avoiding
effect reversal after marginalization. {\it J. R. Statist. Soc. B}, {\bf 65}, 937-941.

\item Cox, D. R. (1972). Regression models and life-tables. {\it J. R. Statist. Soc. B}, {\bf 34}, 187-220.

 \item   Di Serio, C., Rinott, Y. and Scarsini, M. (2009).
Simpson's paradox in survival models . {\it Scand. J.
 Statist.}, {\bf 36}, 463-480.

\item Geng, Z. and Asano, C. (1993). Strong collapsibility of association measures in linear models.
{\it J. R. Statist. Soc. B}, {\bf{55}}, 741-747.

\item Guo, J. H. and Geng, Z. (1995). Collapsibility of logistic regression coefficients.
{\it J. R. Statist. Soc. B}, {\bf 57}, 263-267.

\item   Lindley, D.V. and Novick, M.R. (1981). On the role of
exchangeability in inference. {\it Ann. Statist.}, {\bf 9},
45-58.


\item Ma, Z., Xie, X. and Geng, Z. (2006). Collapsibility of
distribution dependence. {\it J. R. Statist. Soc. B}, {\bf 68},
127-133.


\item Nelson, W. B. (2004). {\it Accelerated Testing: Statistical Models, Test Plans and Data Analysis}.
John Wiley and Sons, New Jersey.


 \item Pearl, J. (1995). Causal diagrams for empirical research. {\it
Biometrika}, {\bf 82}, 110-125.

Pettitt, A. N. (1984). Proportional odds model for survival data
and estimates using ranks. {\it Appl. Statist.}, {\bf 33},
169-175.

\item  Samuels, M.L. (1993). Simpson's paradox and related
phenomena. {\it J.  Amer. Statist. Assoc.}, {\bf
88}, 81-88.

Scarsini, M. and  Spizzichino, F. (1999). Simpson-type paradoxes, dependence, and ageing. {\it J.
Appl. Probab.}, {\bf 36},  119-131.

 \item  Simpson, E.H. (1951). The interpretation of
interaction in contingency tables. {\it J. R. Statist. Soc. B}, {\bf 13}, 238-241.

 \item   Schield, M. (1999). Simpson's paradox and
Cornfield's conditions. {\it Proceedings of the ASA-JSM Section of
Statistical Education, ASA}, 106-111.


\item   Yule, G.U. (1903). Notes on the theory of association of
attributes. {\it Biometrika}, {\bf 2}, 121-134.

\item Vellaisamy, P. and Vijay, V. (2007). Some collapsibility results for n-dimensional contingency tables.
{\it Ann. Inst. Statist. Math.}, {\bf 59}, 557-576.

\item Vellaisamy, P. and Vijay, V. (2008). Collapsibility of
regression coefficients and its extensions. {\it J. Statist.
Plann. Inference}, {\bf 138}, 982-994.

\item Vellaisamy, P. and Vijay, V. (2010). Collapsibility of
contingency tables based on conditional models. {\it J. Statist.
Plann. Inference}, {\bf 140}, 1243-1255.

\item Vellaisamy, P. (2011).  Average collapsibility of distribution
dependence and quantile regression coefficients. To appear in {\it
Scand. J. Statist.}

\item Wermuth, N. (1987). Parametric collapsibility and the lack of moderating effects in contingency tables with a dichotomous response variable. {\it J. R. Statist. Soc. B}, {\bf 49}, 353-364.

\item Wermuth, N. (1989). Moderating effects of subgroups in linear models.
{\it Biometrika}, {\bf 76}, 81-92.

\item Whittemore, A. S. (1978). Collapsibility of multidimensional contingency tables.
{\it J. R. Statist. Soc. B}, {\bf{40}}, 328-340.

\item Whittaker, J. (1990). {\it Graphical Models in Applied Multivariate Statistics}, John Wiley and Sons, New York.

\item   Yule, G.U. (1903). Notes on the theory of association of
attributes. {\it Biometrika}, {\bf 2}, 121-134.

\end{namelist}

\end{document}